\documentclass[a4paper,12pt,reqno]{amsart}
\usepackage{amsmath,amsfonts,amssymb,amsthm,enumerate,multicol}
\usepackage{tikz,graphicx}
\usepackage{stackengine}
\usepackage{hyperref}
\usepackage{caption,enumitem,wasysym}
\usepackage{amsmath,amssymb}
\usepackage{cleveref}
\usepackage{mathrsfs}
\usepackage{physics}
\usepackage{subcaption}
\usepackage{float}
\usepackage{epstopdf} 
\usepackage{graphics,graphicx}
\usepackage{stackengine}
\usepackage{bigints}
\usepackage{geometry}
\usepackage{epsfig}
\usepackage{caption,subcaption}
\usepackage[utf8]{inputenc}
\usepackage{lmodern}
\usepackage[T1]{fontenc}
\usepackage{mathtools}
\usepackage{array}
\usepackage{booktabs}
\usepackage{multirow}
\usepackage{lastpage}
\usepackage{fancyhdr}
\usepackage{pifont}
\usepackage{geometry}
\usepackage{mathtools}
\usepackage{ragged2e}
\usepackage[super]{nth}
\usepackage{float}
\usepackage{float,wrapfig}
\usepackage[labelsep=space]{caption}
\usepackage[font=footnotesize]{caption}
\usepackage{xcolor}
\newtheorem{lemma}{Lemma}
\newtheorem{theorem}{Theorem}[section]
\newtheorem{definition}{Definition}[section]
\newtheorem{example}{Example}[section]
\newtheorem{corollary}{Corollary}[section]

\newtheorem{remark}{Remark}[section]
\newtheorem{note}{Note}[section]

\numberwithin{equation}{section}
\title{Numerical approximation of the Hilfer fractional derivative using orthogonal Quasi-Polynomials}
\author{Rabab Beirkdar\textsuperscript{1}}
\address{\textsuperscript{1} Department of Mathematics\\ Indian Institute of Technology, Roorkee-247667, Uttarakhand, India}
\email{rabab\_ab@ma.iitr.ac.in, rababbeirkdar@gmail.com}
\author[A. Swaminathan]{A. Swaminathan,\textsuperscript{2}}
\address{\textsuperscript{2}Department of Mathematics\\ Indian Institute of Technology, Roorkee-247667, Uttarakhand, India}
\email{a.swaminathan@ma.iitr.ac.in, mathswami@gmail.com}
\allowdisplaybreaks
\bigskip

\begin{document}
\begin{abstract}
A numerical approach based on a composite Gaussian quadrature for the Hilfer fractional derivative is presented. This approach employs quasi-polynomials that are orthogonal with respect to the weight function on the interval $(0,1)$. As special cases, the Riemann--Liouville and Caputo fractional derivatives are obtained. Error bounds for the approximation are derived. To confirm the convergence rates as well as the high efficiency of the new approach, numerical examples are provided.
\end{abstract}
	\maketitle
	\markboth{Rabab Beirkdar and A. Swaminathan}{}
	\section{introduction and preliminaries}\label{introdu}
	The field of fractional calculus extends the concepts of differentiation and integration to arbitrary real orders, providing a refined mathematical framework for systems with memory, non-locality, and anomalous diffusion. Two foundational operators dominate the literature, namely, the Riemann--Liouville fractional derivative \cite{I. Podlubny_1999_Fractional differential equations}, defined as  
	\begin{align*}
	{}^{RL}D^{\alpha}_{a+}f(t)=\frac{1}{\Gamma(n-\alpha)}\frac{d^{n}}{dt^{n}}\int_{a }^{t}(t-\tau)^{n-\alpha-1}f(\tau)\,d\tau, \quad t>0,
\end{align*}
	and the Caputo fractional derivative \cite{I. Podlubny_1999_Fractional differential equations},
	\begin{align*}
	{}^{C}D^{\alpha}_{a+}f(t)=\frac{1}{\Gamma(n-\alpha)}\int_{a}^{t}(t-\tau)^{n-\alpha-1}f^{(n)}(\tau)\,d\tau,  \quad t>0,
	\end{align*}
	where $n-1<\alpha<n$, $n\in \mathbb{N}$, which differ in the order of differentiation and integration.
    To connect these definitions, Hilfer proposed a generalized fractional derivative of a certain orders \(\alpha\) and type \(\beta\in[0,1]\) \cite{R. Hilfer_2000_Applications of fractional calculus in physics}
	\begin{align} \label{lift and right hilfer derivative}
		{}^{H}D^{\alpha,\beta}_{a\pm}f(t)=\pm I_{a\pm}^{\beta(n-\alpha)}\,D^{n}\,I_{a\pm}^{(1-\beta)(n-\alpha)}f(t),
\end{align}
	which reduces to the Riemann--Liouville derivative when \(\beta=0\) and to the Caputo derivative when \(\beta=1\).

	Gaussian quadrature rules for the numerical computation of the   Riemann--Liouville fractional integral were first introduced in \cite{Rapaic_2014_ novel class of fractionally orthogonal quasi-polynomials}, using a class of quasi-polynomials orthogonal with respect to the fractional integration operator. This framework was extended by \cite{G.V. Milovanovic_2021_Some orthogonal polynomials on}, who established the connection between these quasi-polynomials and classical orthogonal polynomials on finite intervals, enabling the use of stable algorithms for the computation of quadrature nodes and weights. The present paper extends this approach to the numerical approximation of the Hilfer fractional derivative \eqref{lift and right hilfer derivative}, which generalizes 
   both the Riemann--Liouville and Caputo derivatives.

 Functions arising naturally in fractional calculus, such as solutions to 
fractional differential equations, often have the form 
of a fractional power series \cite{Rapaic_2014_ novel class of fractionally orthogonal quasi-polynomials}
\begin{align*} 
	f(t) = \sum_{k=0}^{\infty} a_k\, t^{ks}, 
	\quad s > 0.
\end{align*}
Classical orthogonal polynomials are given poor approximations for such functions exhibiting fractional power behavior near $t = 0$. This motivates the use of quasi-polynomials or $s$-polynomials, defined as polynomials in $t^s$ instead of $t$
\begin{align} \label{s-polynomial}
	P_{n,s}^{(\alpha)}(t)
	=\sum_{\nu=0}^{n}c_{n,\nu,s}^{(\alpha)}t^{\nu s}.
\end{align}
Note that when $s=1$, \eqref{s-polynomial} reduces to an ordinary polynomial. For $s = \frac{1}{2}$, the monomials are $\{1, \sqrt{t}, t, t^{3/2}, \ldots\}$, which naturally match the singular behavior of many 
fractional calculus functions.
For numerical evaluation of the left Riemann--Liouville 
fractional integral of order $\alpha > 0$
\begin{align*}
	{}_{0}I_{t}^{\alpha}f(t)
	=\frac{1}{\Gamma(\alpha)}
	\int_{0}^{t}f(\tau)(t-\tau)^{\alpha-1}\,d\tau,
\end{align*}
the quasi-polynomials \eqref{s-polynomial} are required to be orthogonal with respect to the fractional integration operator \cite{Rapaic_2014_ novel class of fractionally orthogonal quasi-polynomials}
\begin{align}	\label{fractional orthogonal}
	{}_{0}I_{1}^{\alpha}\!\left(
	P_{n,s}^{(\alpha)}P_{m,s}^{(\alpha)}
	\right)=0,\quad n\neq m.
\end{align}
As shown in \cite{G.V. Milovanovic_2021_Some orthogonal polynomials on}, this fractional orthogonality condition is equivalent to classical orthogonality with respect to the weight function
\begin{align} \label{weight of two parametrs}
	w(x;\alpha,s) 
	= \frac{(1 - x^{1/s})^{\alpha-1}}
	{\Gamma(\alpha)\,s\, x^{1-1/s}}, 
	\quad \alpha, s > 0, \quad x \in (0,1),
\end{align}
under the substitution $x = t^s$, which transforms \eqref{fractional orthogonal} into
\begin{align*} 	
	\int_0^1 P^{(\alpha)}_n(x)\, 
	P^{(\alpha)}_m(x)\, 
	w(x;\alpha,s)\, dx = 0, 
	\quad n \neq m.
\end{align*}
Since the weight function $w(x;\alpha,s)$ defined in \eqref{weight of two parametrs} is non-negative and integrable on $(0,1)$ for all $\alpha, s > 0$, all the moments
\begin{align*}
	\mu_k = \int_0^1 x^k w(x;\alpha,s)\, dx 
	= \frac{\Gamma(s k + 1)}{\Gamma(\alpha + s k + 1)},
	\quad k = 0, 1, 2, \ldots,
\end{align*}
exist and are finite, the existence of a unique family of monic orthogonal quasi-polynomials satisfying \eqref{fractional orthogonal} follows from \cite{Chihara_1978_ An introduction to orthogonal polynomials}.
	
	The following theorem establishes the existence and recurrence properties of these polynomials.
	\begin{theorem}\cite{Rapaic_2014_ novel class of fractionally orthogonal quasi-polynomials} There exists a family of $s$-polynomials \(P_{n,s}^{(\alpha)}(t)= P_{n}^{(\alpha)}(x)\), $x=t^{s}$, satisfying the orthogonality condition. These quasi-polynomials satisfy the three-term recurrence relation
		\begin{align*}
	P_{k+1}^{(\alpha)}(x)=\Big{(}x-A_{k}^{(\alpha)}\Big{)}\,P_{k}^{(\alpha)}(x)-B_{k}^{(\alpha)}P_{k-1}^{(\alpha)}(x),\quad P_{0}^{(\alpha)}(x)=1,
		\end{align*}
		where the recurrence coefficients are given by
	\begin{align*}
	A_{k}^{(\alpha)}=\frac{\langle xP_{k}^{(\alpha)},P_{k}^{(\alpha)}\rangle_{w}}{\langle P_{k}^{(\alpha)},P_{k}^{(\alpha)}\rangle_{w}},\quad
   B_{k}^{(\alpha)}=\frac{\langle P_{k}^{(\alpha)},P_{k}^{(\alpha)}\rangle_{w}}{\langle P_{k-1}^{(\alpha)},P_{k-1}^{(\alpha)}\rangle_{w}}.
	\end{align*}
		Here, the inner product \(\langle\cdot,\cdot\rangle_{w}\) is defined by
		\begin{align*}
			\langle p,q\rangle_{w}=\int_{0}^{1}p(x)q(x)w(x;\alpha,s)\,dx,
		\end{align*}
		with the weight function
		\begin{align*}
	w(x;\alpha,s)=\frac{1}{\Gamma(\alpha)}\frac{(1-x^{1/s})^{\alpha-1}}{s x^{1-1/s}},\quad\alpha,s>0.
		\end{align*}
\end{theorem}
	\begin{remark} \cite{Rapaic_2014_ novel class of fractionally orthogonal quasi-polynomials}
		For \(s=1\), the first few monic orthogonal polynomials are
		\begin{align*}
			P_{1,1}^{(\alpha)}(t)&=t-\frac{1}{\alpha+1}, \quad
			P_{2,1}^{(\alpha)}(t)=t^{2}-\frac{4t}{\alpha+3}+\frac{2}{(\alpha+2)(\alpha+3)},\\
			P_{3,1}^{(\alpha)}(t)&=t^{3}-\frac{9t^{2}}{\alpha+5}+\frac{18t}{(\alpha+4)(\alpha+5)}-\frac{6}{(\alpha+3)(\alpha+4)(\alpha+5)},\dots
		\end{align*}
		For \(\alpha=1\), the first few quasi-polynomials are
			\begin{align*}
			P_{1,s}^{(1)}(x)&=x-\frac{1}{s+1}, \quad
			P_{2,s}^{(1)}(x)=x^{2}-\frac{2(s+1)x}{3s+1}+\frac{s+1}{(2s+1)(3s+1)},\\
			P_{3,s}^{(1)}(x)&=x^{3}-\frac{3(2s+1)x^{2}}{5s+1}+\frac{3(s+1)(2s+1)x}{(4s+1)(5s+1)}-\frac{(s+1)(2s+1)}{(3s+1)(4s+1)(5s+1)},\dots
		\end{align*} and so forth.
	\end{remark}
	The following theorem provides the Gaussian quadrature rule associated with the weight function \eqref{weight of two parametrs}.
		\begin{theorem} \label{Gauss quadrature rule for approx. the fractional integral}\cite{G.V. Milovanovic_2021_Some orthogonal polynomials on}
		Let $w(x)\equiv w(x;\alpha,s)=(s\Gamma(\alpha))^{-1}x^{1/s-1}(1-x^{1/s})^{\alpha-1}$,$s>0$, with moments
		\begin{align*}
		\mu_{k}=\int_{0}^{1}x^{k}w(x;\alpha,s)\,dx=\frac{\Gamma(s k+1)}{\Gamma(\alpha+s k+1)},\quad k=0,1,\dots .
		\end{align*}
		For each \(n\in\mathbb{N}\), there exists a unique Gauss--Christoffel quadrature rule of maximal degree \(2n-1\),
		\begin{align*}
		\int_{0}^{1}\varphi(x)w(x;\alpha,s)\,dx&=\sum_{k=1}^{n}A_{n,k}(\alpha,s)\,\varphi\big(\xi_{n,k}(\alpha,s)\big)+R_{n}(\varphi;w)\\
		&=Q_n(\varphi;w)+R_n(\varphi;w),
		\end{align*}
		where the nodes $\xi_k=\xi_{n,k}(\alpha,s)$, $k=1,...,n$ are the eigenvalues of the Jacobi matrix and \(A_k=A_{n,k}(\alpha,s)\), $k=1,...,n$ are the corresponding Christoffel numbers. The remainder satisfies \(R_{n}(x^{k};w)=0\) for \(k=0,1,\dots,2n-1\).
		Consequently, for the left Riemann--Liouville fractional integral after we take the transformation \(\tau = t x^{1/s}\),$s>0$ it can approximated by the weighted Gaussian quadrature sum
		\begin{align*}
		{}_{0}I_{t}^{\alpha}f(t)&=t^{\alpha}\sum_{k=1}^{n}A_{n,k}(w)f\big(t\xi_{n,k}^{1/s}\big)+R_{n}(f;w)\\
		&=Q_n(f(t,\alpha,s);w)+R_n(f;w),
		\end{align*}
		for each $f\in C[a,b]$, the sequences of quadrature sums ${Q_n(f(t,\alpha,s);w)}$ converge to ${}_{0}I_{t}^{\alpha}f(t)$ and the rate of convergence is determined by the properties of the function $f$.
	   \end{theorem}
	This quadrature framework provides the foundation for the numerical approximation of fractional integrals, which we extend to the Hilfer derivative in the next section.
		The rest of the paper is organized as follows. 
	Section \ref{numerical} extends the approach to the Hilfer fractional derivative and provides a stepwise numerical procedure. Section \ref{numerical results} demonstrates the method's accuracy through numerical examples.
\section{Main Results: Approximation Formula, Special Cases, and Error Estimates}\label{numerical}
In this section, we represent the Gaussian quadrature framework from Section \ref{introdu} to approximate the Hilfer fractional derivative, providing a clear three-stage process, approximating the inner fractional integral, differentiating, and then approximating the outer integral.\\
For numerical computation, we consider the case $a = 0$ in \eqref{lift and right hilfer derivative}. 
Let $f \in C[0,t]$ with $t > 0$, 
$\alpha \in (0,1)$, $\beta \in [0,1]$, 
and $\gamma = (1-\alpha)(1-\beta)$. 
Under these conditions, the left Hilfer fractional 
derivative has the integral representation
\begin{align} \label{left hilfer integral}
	D_{0+}^{\alpha,\beta} f(t) 
	= \frac{1}{\Gamma(\beta(1-\alpha))\Gamma(\gamma)} 
	\int_0^t (t-\tau)^{\beta(1-\alpha)-1} 
	\left[ \frac{d}{d\tau} \int_0^\tau 
	(\tau-s)^{\gamma-1} f(s)\, ds \right] d\tau.
\end{align}


\begin{theorem} \label{HilferApprox}
Under the conditions \eqref{left hilfer integral}, for any positive integers \(n, m\), the derivative can be approximated by
\begin{align*}
		D_{0+}^{\alpha,\beta}f(t) \approx t^{\beta(1-\alpha)} \sum_{j=1}^{m} B_j \, g\!\left( t \xi_j^{1/s} \right),
\end{align*}
where
	\begin{align*}
		g(\tau) = \sum_{k=1}^{n} A_k \left[ \gamma \tau^{\gamma-1} f\!\left( \tau \xi_k^{1/s} \right) + \tau^{\gamma} f'\!\left( \tau \xi_k^{1/s} \right) \cdot \xi_k^{1/s} \right],
	\end{align*}
	and the following holds.
	\begin{itemize}
		\item \(\xi_k, A_k\) are the nodes and weights of the Gaussian quadrature rule with respect to the weight function
	\begin{align*}
	w(x;\gamma,s) = \frac{x^{\frac{1}{s}-1}}{s \,\Gamma(\gamma)} (1 - x^{1/s})^{\gamma-1}, \quad x \in [0,1],
	\end{align*}
		corresponding to the inner fractional integral \(I^\gamma\).
		\item \(\xi_j, B_j\) are the nodes and weights for the weight function
	\begin{align*}
	w(x;\beta(1-\alpha),s) = \frac{x^{\frac{1}{s}-1}}{s \,\Gamma(\beta(1-\alpha))} (1 - x^{1/s})^{\beta(1-\alpha)-1},
	\end{align*}
	corresponding to the outer fractional integral \(I^{\beta(1-\alpha)}\).
	\end{itemize}
\end{theorem}

\begin{proof}
For the inner fractional integral,  
	we have,
\begin{align*}
		I^\gamma f(t) = \frac{1}{\Gamma(\gamma)} \int_0^t f(\tau) (t-\tau)^{\gamma-1} d\tau.
\end{align*}
	Using the substitution \(\tau = t x^{1/s}\), we obtain,
\begin{align*}
		I^\gamma f(t) = \frac{t^\gamma}{s \Gamma(\gamma)} \int_0^1 f\!\left( t x^{1/s} \right) x^{\frac{1}{s}-1} (1 - x^{1/s})^{\gamma-1} dx.
\end{align*}
  Now the weight function is
  \begin{align*}
  	w(x;\gamma, s) = \frac{x^{\frac{1}{s}-1}}{s \Gamma(\gamma)} (1-x^{1/s})^{\gamma-1}.
  \end{align*}
  \begin{align*}
  	\Longrightarrow I^{\gamma}f(t) = t^{\gamma} \int_{0}^{1} f(tx^{1/s}) w(x; \gamma, s) \, dx.
  \end{align*}
By using Theorem \ref{Gauss quadrature rule for approx. the fractional integral}, this fractional integral can be approximated by the weighted Gaussian quadrature sum as
	\begin{align} \label{integralaprox.}
		I^{\gamma}f(t) = t^{\gamma} \sum_{k=1}^{n} A_k f(t \xi_{k}^{1/s}),
	\end{align}
	where the nodes \(\xi_{k}\) and the weights \(A_k\), \(k=1,\ldots,n\), can be calculated by using the Mathematica\textsuperscript{\copyright}software in specific, we used the package \texttt{OrthogonalPolynomials} \cite{Milovanovic_2004_ OrthogonalPolynomials,Milovanovic_2012_orthogonal polynomials and corresponding quadratures}, where this package can obtain the nodes and weights in a procedure with a sufficiently high precision.\\
	Now, by taking the first derivative for the both sides of the approximation \eqref{integralaprox.} with respect to \(t\), we get
	\begin{align*}
		\frac{d}{dt} I^\gamma f(t) &\approx \frac{d}{dt} \left[ t^\gamma \sum_{k=1}^{n} A_k f\!\left( t \xi_k^{1/s} \right) \right] 
		&= \sum_{k=1}^{n} A_k \left[ \gamma t^{\gamma-1} f\!\left( t \xi_k^{1/s} \right) + t^\gamma f'\!\left( t \xi_k^{1/s} \right) \cdot \xi_k^{1/s} \right] = g(t),
	\end{align*}
	as given in the hypothesis. For outer fractional integral, the Hilfer derivative becomes
\begin{align*}
		D_{0+}^{\alpha,\beta} f(t) = I^{\beta(1-\alpha)} g(t) = \frac{1}{\Gamma(\beta(1-\alpha))} \int_0^t (t-\tau)^{\beta(1-\alpha)-1} g(\tau) \, d\tau.
\end{align*}
	Applying the same substitution \(\tau = t x^{1/s}\), we obtain
	\begin{align*}
	I^{\beta(1-\alpha)} g(t) = t^{\beta(1-\alpha)} \int_0^1 g\!\left( t x^{1/s} \right) w(x;\beta(1-\alpha),s) \, dx.
	\end{align*}
	Using Gaussian quadrature again with nodes \(\xi_j\) and weights \(B_j\) for the weight function \(w(x;\beta(1-\alpha),s)\), we approximate
\begin{align*}
		I^{\beta(1-\alpha)} g(t) \approx t^{\beta(1-\alpha)} \sum_{j=1}^{m} B_j \, g\!\left( t \xi_j^{1/s} \right).
\end{align*}
Thus we get the final approximation,  
\begin{align*}
D_{0+}^{\alpha,\beta} f(t) \approx t^{\beta(1-\alpha)} \sum_{j=1}^{m} B_j \sum_{k=1}^{n} A_k \left[ \gamma (t\xi_j^{1/s})^{\gamma-1} f\!\left( t\xi_j^{1/s} \xi_k^{1/s} \right) + (t\xi_j^{1/s})^\gamma f'\!\left( t\xi_j^{1/s} \xi_k^{1/s} \right) \xi_k^{1/s} \right],
\end{align*}
which completes the proof.  
\end{proof}
\begin{corollary} \label{RiemannApprox.}
	Let \(\alpha \in (0,1)\) and \(f \in C^{1}[0,t]\). The left-sided Riemann-Liouville fractional derivative corresponds to the Hilfer derivative with \(\beta = 0\) is given by
	\begin{align*}
	D_{0+}^{\alpha,0}f(t) = {}^{RL}D_{0+}^{\alpha}f(t) = \frac{d}{dt} I_{0+}^{1-\alpha} f(t).
	\end{align*}
	Setting \(\beta = 0\) and \(\gamma = 1-\alpha\) in Theorem \ref{HilferApprox}, the numerical approximation simplifies to:
	\begin{align*}
	{}^{RL}D_{0+}^{\alpha}f(t) \approx \sum_{k=1}^{n} A_k \left[ \gamma t^{\gamma-1} f\!\left( t \xi_k^{1/s} \right) + t^\gamma f'\!\left( t \xi_k^{1/s} \right) \cdot \xi_k^{1/s} \right],
	\end{align*}
	and the nodes and weights \(\xi_k, A_k\) correspond to the weight function \(w(x;1-\alpha,s)\).
\end{corollary}

\begin{corollary} \label{CaputoApprox.}
	Let \(\alpha \in (0,1)\) and \(f \in C^{1}[0,t]\). The left-sided Caputo fractional derivative corresponds to the Hilfer derivative with \(\beta = 1\) is given by
	\begin{align*}
	D_{0+}^{\alpha,1}f(t) = {}^CD_{0+}^{\alpha}f(t) = I_{0+}^{1-\alpha} f'(t).
	\end{align*}
	Setting \(\beta = 1\) and \(\gamma = 0\) in Theorem \ref{HilferApprox}, the numerical approximation becomes:
	\begin{align*}
	{}^CD_{0+}^{\alpha}f(t) \approx t^{1-\alpha} \sum_{j=1}^{m} B_j \, f'\!\left( t \xi_j^{1/s} \right),
	\end{align*}
where the nodes \(\xi_j\) and weights \(B_j\) are corresponding to the weight function \(w(x;1-\alpha,s)\).
\end{corollary}
It may be noted that corollaries \ref{RiemannApprox.}, \ref{CaputoApprox.}, which are special cases of Theorem \ref{HilferApprox} is not explicitly available in the literature. Hence these results for the Riemann-Liouville operator and Caputo operator, respectively, attract special attention.

Next, we will find the error estimate for the Hilfer fractional derivative by the proposed method.
We assume that
\begin{equation} \label{assumption of f}
f(\tau) = \tau^{p-1} f_1(\tau), \quad f_1 \in C^1[0,t], \quad p>1.
\end{equation}

\begin{theorem}\cite{Milovanovic_2008_Interpolation processes}\label{error bound for general weight}
 For a general weight function $w$ and $f \in C[-1,1]$, we have:
\begin{align*}
|R_n(f)_w| \le 2 \|w\|_1 \, E_{2n-1}(f)_\infty,
\end{align*}
where
\begin{align*}
\|w\|_1 = \int_{-1}^{1} w(x)\,dx.
\end{align*}
\end{theorem}

\begin{remark}
	Theorem~\ref{error bound for general weight} holds for the weight function $w(x;\alpha,\beta)$ defined in \eqref{weight of two parametrs} on the interval $(0,1)$, since $w(x;\alpha,\beta)$ is non-negative 
	and integrable on $(0,1)$ for all $\alpha, \beta > 0$, with 
\begin{align*}
	\|w\|_1 = \int_{0}^{1} w(x)\,dx.
\end{align*}
\end{remark}

\begin{lemma} \label{UPA} Let $f$ be as the assumption \eqref{assumption of f} and define the transformed function:
	\begin{align*}
g(x) = f(t x^{1/s}), \quad 0 < s \leq 1, \quad x \in [0,1], \quad t > 0.
	\end{align*}
Let $E_{2n-1}(g)$ denote the best uniform polynomial approximation of degree $2n-1$ on $[0,1]$. Then there exists a constant $C>0$ independent of $f_1$ and t such that:
\begin{align*}
E_{2n-1}(g) \le Ct^{p-1}\|f_1\|_{C^1} (2n-1)^{-\sigma}, \quad \sigma = \min\left(1, \frac{p-1}{s}\right).
\end{align*}
\end{lemma}
\begin{proof}
Let 
\begin{align*}
g(x) = t^{p-1} x^{\frac{p-1}{s}} f_1(t x^{1/s}).
\end{align*}
Since $f_1 \in C^1[0,t]$, the function $f_1(t x^{1/s})$ is smooth in $[0,1]$. Hence, The factor $x^\mu$ controls the regularity of $g(x)$. Where $\mu = \frac{p-1}{s}.$
That is,
$g(x) \sim x^\mu$ and \\ $\|f_1\|_{C^1} = \max\left(\|f_1\|_{\infty},\ \|f_1'\|_{\infty}\right)$.
\begin{enumerate}
	\item If $\mu \geq 1$, the function $g(x) \in C^1[0,1]$.  
	By Jackson's theorem \cite[p.35]{Milovanovic_2008_Interpolation processes}:
	\begin{align*}
	E_m(g) \leq Ct^{p-1}\|f_1\|_{C^1} m^{-1}.
	\end{align*}
	
	\item If $0 < \mu < 1$, the function $g(x)$ is H\"older continuous of order $\mu$.  
	By Jackson's theorem \cite[p.35]{Milovanovic_2008_Interpolation processes}:
	\begin{align*}
	E_m(g) \leq Ct^{p-1}\|f_1\|_{C^1} m^{-\mu}.
	\end{align*}
\end{enumerate}
Combining both cases, we obtain:
\begin{align*}
E_m(g) \leq Ct^{p-1}\|f_1\|_{C^1} m^{-\min(1,\mu)}.
\end{align*}
Setting $m = 2n - 1$ completes the proof.
\end{proof}
\begin{theorem}
 Let $f$ satisfy the above assumptions given in \eqref{assumption of f} with $p > 2-\gamma$, $0 < s \leq 1$, $0<\alpha,\beta<1$ and let $Q_{n,m}(f;t)$ be the proposed quadrature approximation of the Hilfer fractional derivative. Then the error is estimated as
 \begin{align*}
 	\left|{}^{H}D^{\alpha,\beta}_{0^+}f(t) - Q_{n,m}(f;t)\right| 
 	\leq  
 	\frac{Ct^{\beta(1-\alpha)}}{\Gamma(\beta(1-\alpha)+1)}
 	\|f_1\|_{C^1}
 	\left[
 	\frac{(2n-1)^{-\min(1,\mu)}}{\Gamma(\gamma+1)} 
 	+ 
 	(2m-1)^{-\min(1,\mu_{*})}
 	\right],
 \end{align*}
where $\mu = \dfrac{p-1}{s}$, $\gamma = (1-\alpha)(1-\beta)$ and $\mu_{*} = \dfrac{p+\gamma-2}{s}$.
\end{theorem}
\begin{proof}
Using the definition of the Hilfer fractional derivative, we write
\begin{align*}
D^{\alpha,\beta} f(t)
= I^{\beta(1-\alpha)} g(t),
\quad
g(t) = \frac{d}{dt} I^\gamma f(t),
\end{align*}
where $\gamma = (1-\alpha)(1-\beta)$.
Let $g_n(t)$ be the approximation of $g(t)$ obtained from the inner quadrature rule.  
Hence, the error can be written as
\begin{align*}
	D^{\alpha,\beta}f(t) - Q_{n,m}(f;t) 
	= 
	\underbrace{I^{\beta(1-\alpha)}(g - g_n)(t)}_{\text{inner error}} 
	+ 
	\underbrace{\left(I^{\beta(1-\alpha)}g_n(t) - Q_m(g_n;t)\right)}_{\text{outer error}}.
	\end{align*}
We will estimate the bound for each term separately.
Using the definition of the fractional integral, we have
\begin{align*}
\left| I^{\beta(1-\alpha)} (g - g_n)(t) \right|
\le
\frac{1}{\Gamma(\beta(1-\alpha))}
\int_0^t (t-\tau)^{\beta(1-\alpha)-1}
|g(\tau) - g_n(\tau)| \, d\tau.
\end{align*}
Therefore, by the boundedness of the Riemann--Liouville fractional integration operator on $C[0,t]$ \cite[Lemma~3.3, p.~163]{A. Kilbas_2006_ Theory and applications of fractional differential equations}, we get 
\begin{align} \label{g-gn}
	\left|I^{\beta(1-\alpha)}(g - g_n)(t)\right| 
\le
\frac{t^{\beta(1-\alpha)}}{\Gamma(\beta(1-\alpha)+1)}
\|g - g_n\|_\infty.
\end{align}
From the definition of $g(t)$ and its quadrature approximation $g_n(t)$, we have
\begin{align*}
g(t) - g_n(t)
= \gamma t^{\gamma-1} R_n(F) + t^\gamma R_n(F_1),
\end{align*}
where
\begin{align*}
F(x) = f(t x^{1/s}),
\quad
F_1(x) = x^{1/s} f'(t x^{1/s}).
\end{align*}
Using Theorem \ref{error bound for general weight}, we obtain
\begin{align*}
|R_n(F)| \le C \|w\|_1 E_{2n-1}(F)
\quad \mbox{and} \quad
|R_n(F_1)| \le C \|w\|_1 E_{2n-1}(F_1).
\end{align*}
Hence,
\begin{align} \label{bound}
	|g(t) - g_n(t)|
	\le
	C \|w\|_1 \left( \gamma t^{\gamma-1} E_{2n-1}(F)
	+ t^\gamma E_{2n-1}(F_1) \right).
\end{align}
By Lemma \ref{UPA}, both $F$ and $F_1$ satisfy
\begin{align*}
E_{2n-1}(\cdot) \le Ct^{p-1} \|f_1\|_{C^1} (2n-1)^{-\min(1,\mu)},
\end{align*}
where $\mu = \dfrac{p-1}{s}$. By substituting in \eqref{bound}, we get
\begin{align*}
|g(t)-g_n(t)|
\le
C \|w\|_1 \|f_1\|_{C^1}
	\left[
\gamma t^{\gamma-1} \cdot t^{p-1} + t^{\gamma} \cdot t^{p-1}
\right]
(2n-1)^{-\min(1,\mu)}
\end{align*}
\begin{align*}
	= 
	C \|w\|_1\|f_1\|_{C^1}(2n-1)^{-\min(1,\mu)}
	\left[
	\gamma t^{p+\gamma-2} + t^{p+\gamma-1}
	\right].
\end{align*}
Taking the supremum over $t \in [0,T]$, since $p > 2 - \gamma$ then $p + \gamma - 2 > 0$, both powers are bounded at $t = 0$, 
and we obtain
\begin{align*}
	\|g - g_n\|_\infty 
	\leq 
	C  \|w\|_1\|f_1\|_{C^1}(2n-1)^{-\min(1,\mu)}
	\left[
	\gamma T^{p+\gamma-2} + T^{p+\gamma-1}
	\right].
\end{align*}
Since $\|w\|_1 =\dfrac{1}{\Gamma(\gamma+1)}$, 
substituting into \eqref{g-gn}, the inner error satisfies
\begin{align} \label{inner_bound}
	\left|I^{\beta(1-\alpha)}(g - g_n)(t)\right| 
	\leq 
	\frac{Ct^{\beta(1-\alpha)}}{\Gamma(\beta(1-\alpha)+1)\,\Gamma(\gamma+1)}
	\|f_1\|_{C^1}
	(2n-1)^{-\min(1,\mu)}.
\end{align}
We now find the bound for the outer error, which comes from the m-point Gaussian quadrature rule which is applied to approximate
 $I^{\beta(1-\alpha)}g_n(t)$. Using the substitution $\tau = tx^{1/s}$
\begin{align}\label{outer_integral}
	I^{\beta(1-\alpha)}g_n(t) 
	&= \frac{1}{\Gamma(\beta(1-\alpha))} 
	\int_0^t (t-\tau)^{\beta(1-\alpha)-1} g_n(\tau)\, d\tau \notag\\
	&= \frac{t^{\beta(1-\alpha)}}{s\,\Gamma(\beta(1-\alpha))} 
	\int_0^1 x^{1/s-1}
	\left(1 - x^{1/s}\right)^{\beta(1-\alpha)-1} 
	g_n\!\left(tx^{1/s}\right) dx \notag\\
	&= t^{\beta(1-\alpha)} 
	\int_0^1 G_n(x)\, w(x;\beta(1-\alpha),s)\, dx,
\end{align}
where we define
\begin{align*}
	G_n(x) = g_n\!\left(tx^{1/s}\right),\quad x \in [0,1],
\end{align*}
and the weight function is
\begin{align*}
	w(x;\beta(1-\alpha),s) 
	= \frac{x^{1/s - 1}}{s\,\Gamma(\beta(1-\alpha))}
	\left(1 - x^{1/s}\right)^{\beta(1-\alpha)-1},\quad \alpha,\beta,s > 0.
\end{align*}
The $m$-point Gaussian quadrature rule associated with the weight function 
$w(x;\beta(1-\alpha),s)$ approximates the integral in \eqref{outer_integral} as
\begin{align}\label{gauss_outer}
	\int_0^1 G_n(x)\, w(x;\beta(1-\alpha),s)\, dx 
	= \sum_{j=1}^m B_j\, G_n(\xi_j) + R_m(G_n; w_{*}),
\end{align}
where $\xi_j = \xi_{m,j}(\beta(1-\alpha),s)$ are the 
quadrature nodes, $B_j = B_{m,j}(\beta(1-\alpha),s)$ 
are the corresponding Christoffel numbers, and
\begin{align*}
	R_m(G_n; w_{*}) 
	= \int_0^1 G_n(x)\, w_{*}(x)\, dx 
	- \sum_{j=1}^m B_j\, G_n(\xi_j)
\end{align*}
is the remainder term of the quadrature rule.
Consequently, the composite quadrature sum 
$Q_m(g_n;t)$ is given by
\begin{align} \label{eq:Qm_def}
	Q_m(g_n;t) 
	= t^{\beta(1-\alpha)} \sum_{j=1}^m B_j\, G_n(\xi_j).
\end{align}
Subtracting \eqref{eq:Qm_def} from \eqref{outer_integral} 
and using \eqref{gauss_outer}, we obtain
\begin{align}
	I^{\beta(1-\alpha)}g_n(t) - Q_m(g_n;t) 
	&= t^{\beta(1-\alpha)} 
	\int_0^1 G_n(x)\, w_{*}(x)\, dx 
	- t^{\beta(1-\alpha)} \sum_{j=1}^m B_j\, G_n(\xi_j) 
	\notag\\
	&= t^{\beta(1-\alpha)} 
	\left[
	\int_0^1 G_n(x)\, w_{*}(x)\, dx 
	- \sum_{j=1}^m B_j\, G_n(\xi_j)
	\right] \notag\\
	&= t^{\beta(1-\alpha)}\, R_m(G_n; w_{*}).
	\label{eq:outer_diff}
\end{align}
Since $t > 0$ and $\beta(1-\alpha) > 0$, 
we have $t^{\beta(1-\alpha)} > 0$, 
and by taking the absolute value 
of \eqref{eq:outer_diff} gives
\begin{align*}
	\left|I^{\beta(1-\alpha)}g_n(t) - Q_m(g_n;t)\right| 
	= t^{\beta(1-\alpha)}\, 
	\left|R_m(G_n; w_{*})\right|.
\end{align*}
Since $\|w_{*}\|_1 = \frac{1}{\Gamma(\beta(1-\alpha)+1)}$, by Theorem \ref{error bound for general weight}
\begin{align*}
	\left|R_m(G_n)\right| 
	\leq 
	\frac{C}{\Gamma(\beta(1-\alpha)+1)}\,E_{2m-1}(G_n).
\end{align*}
Now we apply Lemma~\ref{UPA} to $G_n(x) = g_n(tx^{1/s})$. 
Since $g_n(t)$ approximates $g(t) = \frac{d}{dt}I^{\gamma}f(t)$, 
and $I^{\gamma}f(t) \sim t^{p+\gamma-1}$ near $t = 0$, 
differentiating gives $g(t) \sim t^{p+\gamma-2}$. 
Therefore substituting $t \to tx^{1/s}$:
\begin{align*}
	G_n(x) = g_n(tx^{1/s}) \sim \left(tx^{1/s}\right)^{p+\gamma-2} 
	\sim x^{\frac{p+\gamma-2}{s}} = x^{\mu_{\mathrm{*}}}.
\end{align*}
By Lemma~\ref{UPA} applied to $G_n$
\begin{align}\label{outer_bound}
	\left|I^{\beta(1-\alpha)}g_n(t) - Q_m(g_n;t)\right| 
	\leq 
	\frac{Ct^{\beta(1-\alpha)}}{\Gamma(\beta(1-\alpha)+1)}
	\|f_1\|_{C^1}
	(2m-1)^{-\min(1,\mu_{*})}.
\end{align}
Adding the inner bound~\eqref{inner_bound} 
and the outer bound~\eqref{outer_bound}:
\begin{align*}
	\left|D^{\alpha,\beta}f(t) - Q_{n,m}(f;t)\right| 
	\leq 
	\frac{Ct^{\beta(1-\alpha)}}{\Gamma(\beta(1-\alpha)+1)}
	\|f_1\|_{C^1}
	\left[
	\frac{(2n-1)^{-\min(1,\mu)}}{\Gamma(\gamma+1)} 
	+ 
	(2m-1)^{-\min(1,\mu_{*})}
	\right].
\end{align*}
This completes the proof.
\end{proof}
\begin{remark}
Since $\mu_{*} < \mu$ always, as $\gamma = (1-\alpha)(1-\beta) < 1$.  $\mu_{*}$ determines the total convergence rate.
The inner error achieves the optimal rate $\mathcal{O}((2n-1)^{-1})$ 
when $s \leq p - 1$, while the outer error achieves the optimal rate 
$\mathcal{O}((2m-1)^{-1})$ when $s \leq p + \gamma - 2$. The assumption $p > 2 - \gamma$ must be satisfied to ensure	$\mu_{*} > 0$ and the validity of the error bound. 
\end{remark}
\begin{remark}
	For the Caputo case $\beta = 1$, the method reduces to approximating 
$I^{1-\alpha}f'(t)$ directly, and performs well only when $f'$ 
is smooth at the origin.
\end{remark}

\section {Numerical Examples: Convergence and Comparison of Fractional Derivative Approximations} \label{numerical results}
\begin{figure}[htb!]
	\centering
	\begin{minipage}[b]{0.48\textwidth}
		\includegraphics[width=\textwidth]{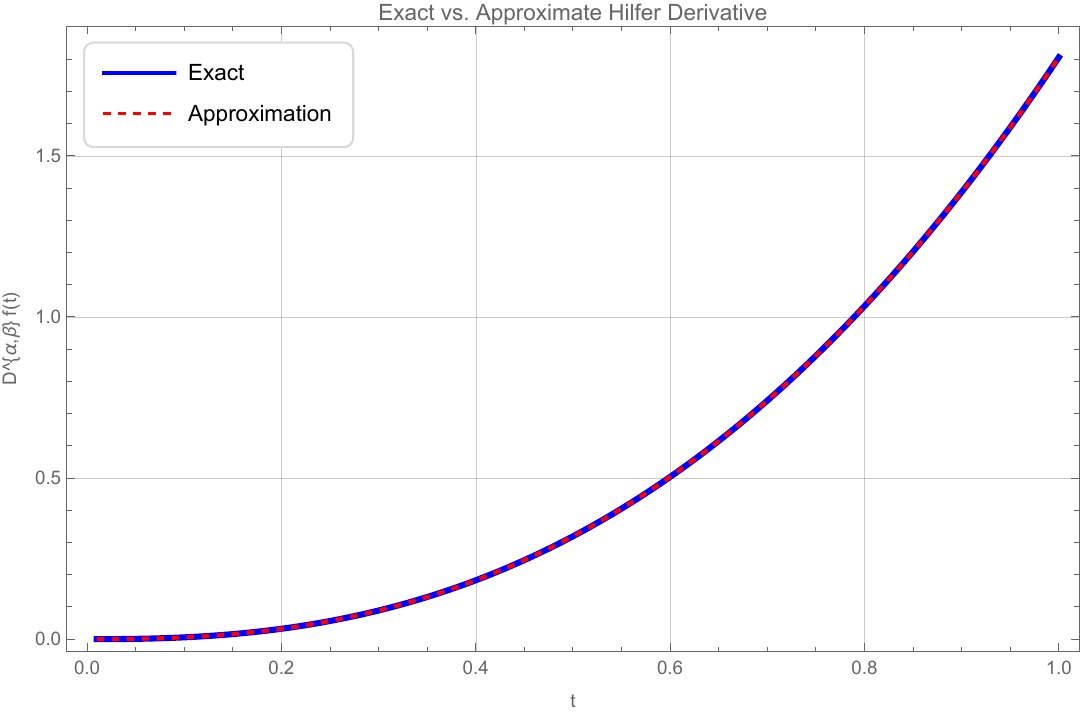}
		\caption{Approximate and Exact Solutions of Ex.\ref{ex1} for $\alpha=1/2$, $\beta=1/3$, $\lambda=3$, $m=n=5$.}
		\label{fig1}
	\end{minipage}
	\hfill
	\begin{minipage}[b]{0.50\textwidth}
		\includegraphics[width=\textwidth]{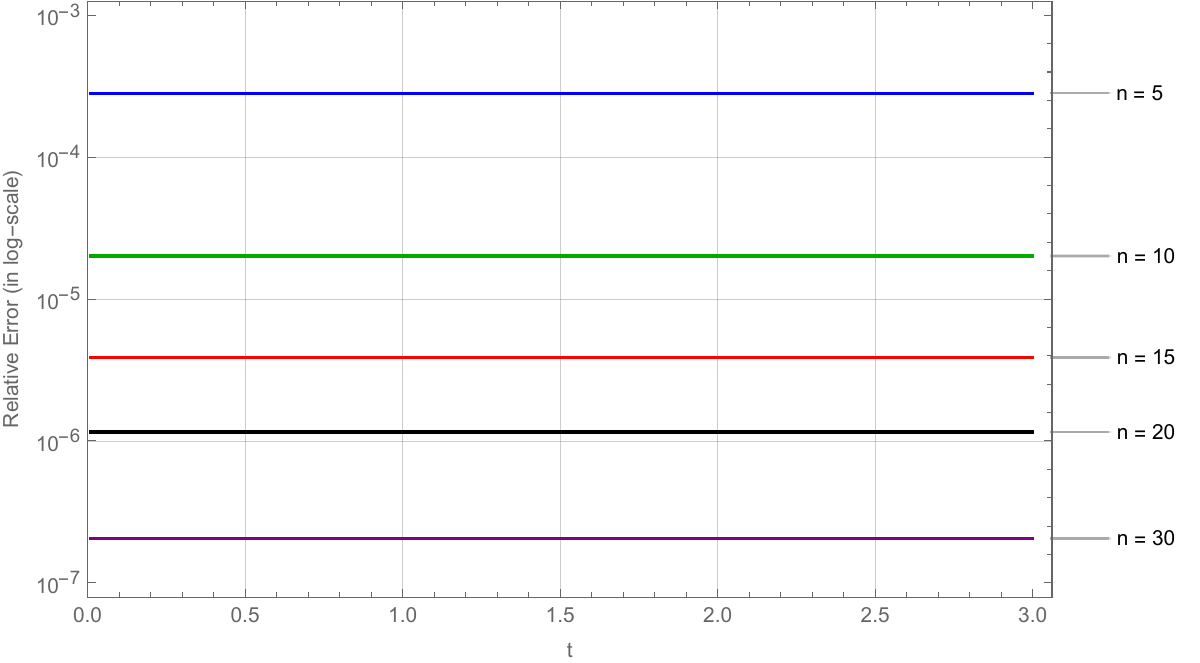}
		\caption{Relative error for different values of n=m for Ex.\ref{ex1} with $\alpha=1/2$, $\lambda=1/2$, $\beta=1/2$(Hilfer), $s=1/3$.}
		\label{fig2}
	\end{minipage}
\end{figure}

\begin{figure}[htb!]
	\centering
	\begin{minipage}[b]{0.48\textwidth}
		\includegraphics[width=\textwidth]{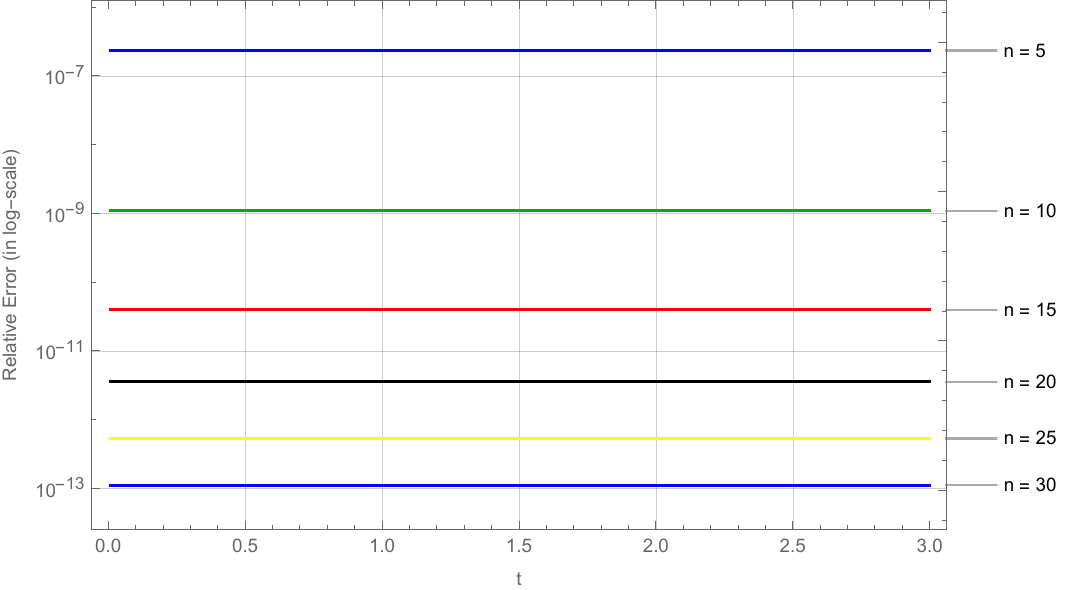}
		\caption{Relative error for different values of n=m for Ex.\ref{ex1} with $\alpha=1/2$, $\lambda=1/2$, $\beta=0$(Riemann), $s=1/3$.}
		\label{fig3}
	\end{minipage}
	\hfill
	\begin{minipage}[b]{0.50\textwidth}
		\includegraphics[width=\textwidth]{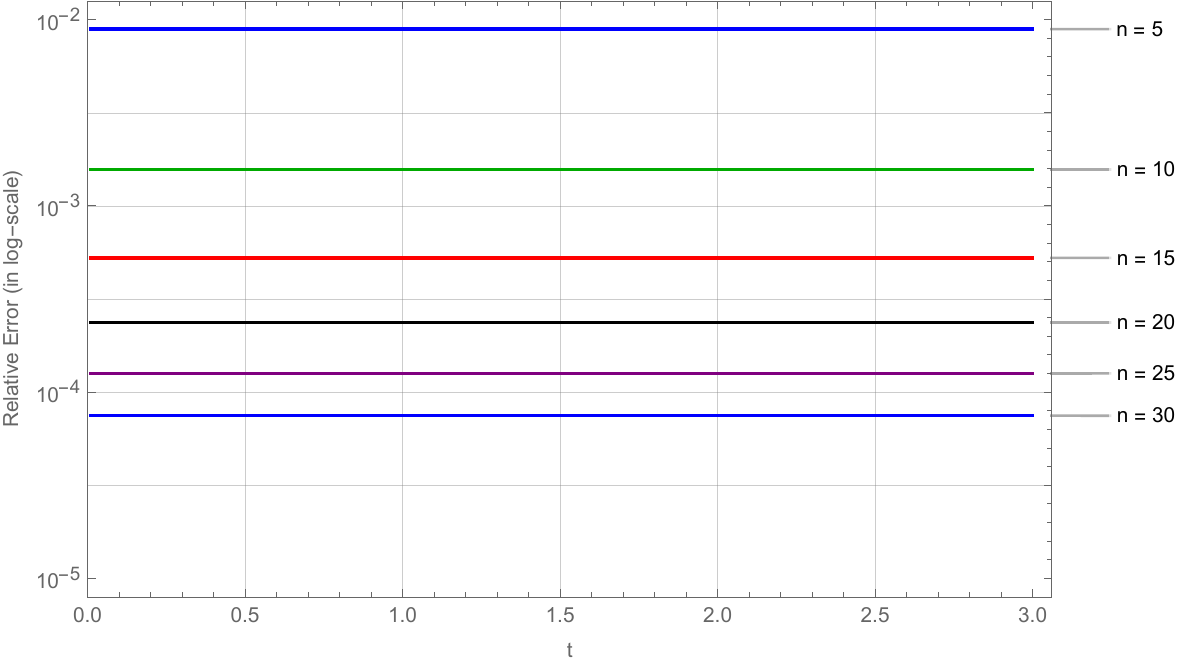}
		\caption{Relative error for different values of n=m for Ex.\ref{ex1} with $\alpha=1/2$, $\lambda=1/2$, $\beta=1$(Caputo), $s=1/3$.}
		\label{fig4}
	\end{minipage}
	\label{fig:combined}
\end{figure}

In this part, we intend to verify the effectiveness of our final proposed formula for approximate Hilfer fractional derivative then compare the results with Riemann-Liouville fractional derivative and Caputo fractional derivative.\\
In all examples, we calculate the relative errors as,\\
\begin{align*}
\mathrm{Err}_{n}^{(\alpha,\beta,s)}f(t) = \left| \frac{Q_{n,m}(f;t) - D_{0+}^{\alpha,\beta}f(t)}{D_{0+}^{\alpha,\beta}f(t)} \right|, \quad (a \leq t \leq b),
\end{align*}
\\
where the numerical approximation derived in Theorem \ref{HilferApprox}  is denoted by the composite Gaussian quadrature sum $Q_{n,m}(f;t)$.
\\
To calculate the exact solution of Hilfer fractional derivative of a function we use the following definition.
\begin{definition} \label{laplace}\cite{R. Hilfer_2000_Fractional time evolution,R. Hilfer_2000_Applications of fractional calculus in physics} 
The Laplace transform of the Hilfer fractional derivative :
\begin{align*}
\mathcal{L}\left\{D_{a+}^{\alpha,\beta} f(x)\right\}(\sigma) = \sigma^{\alpha} \mathcal{L}\{f(x)\}(\sigma) - \sigma^{\beta(\alpha-1)}\left( I_{0+}^{(1-\beta)(1-\alpha)} f \right) (0^+) ,
\end{align*}

where 
\begin{align*}
\left( I_{0+}^{(1-\beta)(1-\alpha)} f \right) (0^+)
\end{align*}
is the Riemann--Liouville fractional integral of order $(1 - \beta)(1 - \alpha)$ evaluated in the limit as $t \to 0^+$.
\end{definition}
	\begin{example} \label{ex1}
	Consider the function $f(t)=t^\lambda$, 
	where the exact Hilfer fractional derivative is given by 
	\begin{align*}
	{}^HD^{\alpha, \beta}_{0} t^{\lambda}=	\frac{\Gamma(1+\lambda)}{\Gamma(1+\lambda-\alpha)} t^{\lambda-\alpha} ,\ \alpha>0,\ \lambda>-1.
	\end{align*}
	
\end{example}
\begin{table}[htbp]
	\centering
	\caption{Relative Errors (RE) for Example \ref{ex1} at $t=0.5$ with $\lambda=3$, $\alpha=1/2, s=1/3$.}
\label{table1}
	\begin{tabular}{cccc}
		\toprule
		$n=m$ & Hilfer ($\beta=1/3$)  & Riemann ($\beta=0$) & Caputo ($\beta=1$)  \\
		\midrule
		2 & 6.75746(-2) & 7.92662(-2) & 2.04748(-2) \\
		3 & 3.06937(-3) & 4.34495(-3) & 1.57359(-4) \\
		4 & 3.65115(-5) & 5.82459(-5) & 0 \\
		5 & 0 & 0 & 0 \\
		\bottomrule
	\end{tabular}
\end{table}
\begin{note}
	In Table \ref{table1} and other tables in the sequel, numbers in parentheses indicate decimal exponents.
\end{note}
In Table \ref{table1}, where we take $\lambda=3$, $\alpha=1/2$ and $s=1/3$, the quadrature process provides excellent approximations for all three fractional derivatives. We can see that for m=n=5, the relative error drops to machine zero, confirming the theoretical exactness of the quadrature rule for a polynomial of degree $3$. Figure \ref{fig1} shows that the exact solution matches the approximation.\\
Now, when we take $\lambda=1/2$, we can see from figures \ref{fig2}, \ref{fig3}, \ref{fig4} for the Hilfer derivative, Riemann derivative and Caputo derivative respectively, relative error almost constant for each t in the interval $(0,3)$, the convergence of the composite Gaussian quadrature sum for the Hilfer fractional derivative is satisfactory and slow for the Caputo derivative. But it is showing the best for the Riemann-Liouville derivative.
\begin{example} \label{ex2}
	Consider the left Hilfer fractional derivative of the function $f(t)=e^t \operatorname{erfc}(\sqrt{t})$, $t>0$ , with $\operatorname{erfc}$ the well-known complementary error-function. By applying the Definition \ref{laplace}, we get the solution as
	\begin{align*}
		{}^HD^{\alpha, \beta}_{0} e^t\operatorname{erfc}(\sqrt{t})= t^{-\alpha} E_{1/2,1-\alpha}(-\sqrt t),
	\end{align*}
	where $E_{\alpha,\beta}(z)$ is the two-parametric Mittag-leffler function defined by
	\begin{align*}
	E_{\alpha,\beta}(z) = \sum_{k=0}^\infty \frac{z^k}{\Gamma(\alpha k + \beta)},\quad \alpha > 0, \quad 
	\beta \in \mathbb{C}, \quad z \in \mathbb{C}.
	\end{align*}
	\end{example}
The function \( f(t) = e^t \operatorname{erfc}(\sqrt{t}) \) has a critical singularity at \( t = 0 \) . The series expansion at the origin is given by
\begin{align*}
f(t) = 1 - \frac{2\sqrt{t}}{\sqrt{\pi}} + t - \frac{4t^{3/2}}{3\sqrt{\pi}} + \frac{t^2}{2} - \frac{8t^{5/2}}{15\sqrt{\pi}} + \mathcal{O}(t^{3}).
\end{align*}

\begin{figure}[htb!]
	\centering
	\begin{minipage}[b]{0.48\textwidth}
		\includegraphics[width=\textwidth]{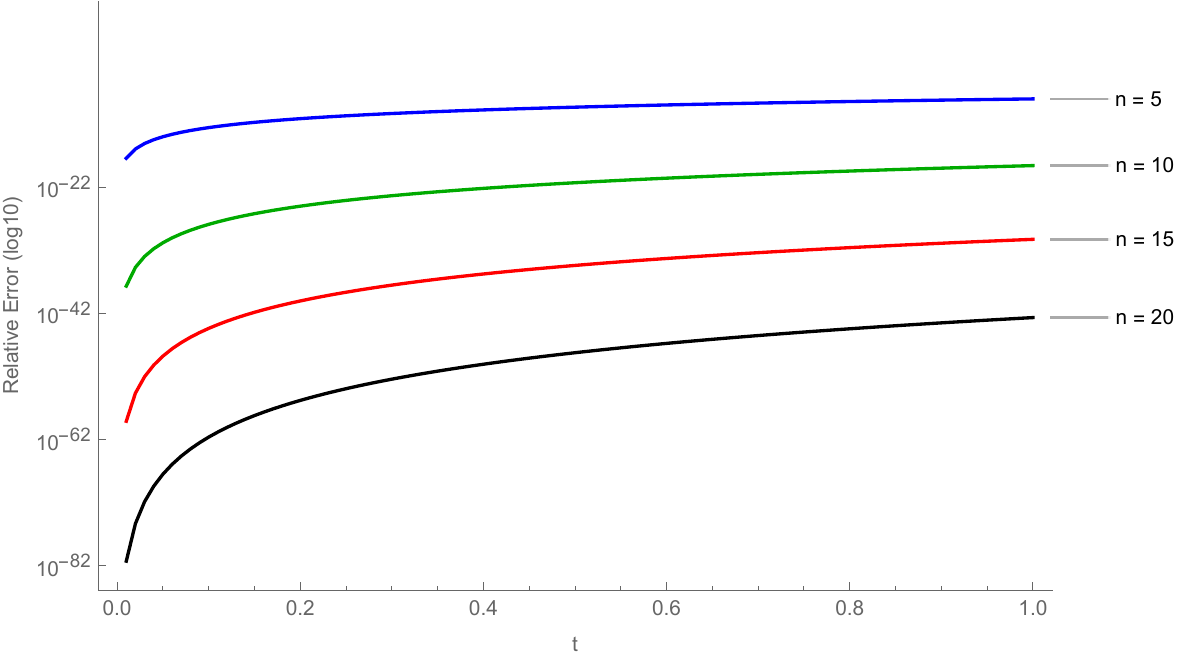}
			\caption{The relative errors for Riemann fractional derivative ($\beta=0$) for Example \ref{ex2} when  $\alpha=1/4$, $s=1/2$ and different values of n.}
		\label{fig5}
	\end{minipage}
	\hfill
	\begin{minipage}[b]{0.48\textwidth}
		\includegraphics[width=\textwidth]{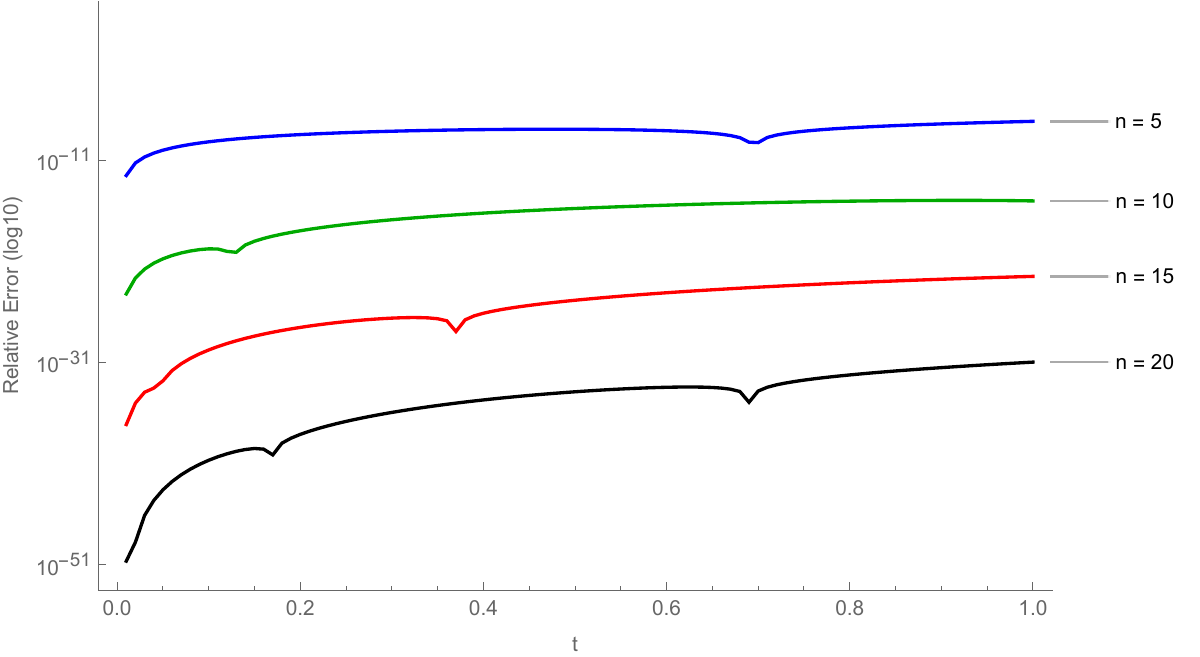}
			\caption{The relative errors for Riemann fractional derivative ($\beta=0$) for Example \ref{ex2} when  $\alpha=1/4$, $s=1/4$ and different values of n.}
		\label{fig6}
	\end{minipage}
\end{figure}

Because of the singularity, convergence will be slow, hence choosing proper $s$ and $\beta$, we can remove this singularity. We can see from Table \ref{table2} and Figure \ref{fig5} for $\alpha=1/4$ and $s=1/2$ for the Hilfer ($\beta=1/4$) and the Caputo fractional derivatives, convergence is slow, but for the Riemann fractional derivative, it is very fast. \\
Next, when we take $s=1/4$, we can see from Table \ref{table3} and Figure \ref{fig6} the convergence showing better for the Hilfer derivative, but for the Caputo derivative it is still slow, and for the Riemann it converges fast when we increase the number of $n=m$.
\begin{table}[h]
	\centering
	\small
	\caption{Relative errors in composite Gaussian quadrature sum obtained $n$-point for fractional derivative Hilfer (H($\beta=1/4$)) and Caputo(C($\beta=1$)) in some selected values of $t \in [0, 1]$ for $\alpha = 1/4$ and $s = 1/2$.} \label{table2}
	\begin{tabular}{llllllll}
		\toprule
		\textbf{$n=m$} & \textbf{FD}  &   \textbf{$t=0.1$} & \textbf{$t=0.2$} & \textbf{$t=0.3$} & \textbf{$t=0.5$} & \textbf{$t=0.8$} & \textbf{$t=1.0$} \\ \midrule
		5  & H($\beta=1/4$)   & 7.52(-3) & 8.66(-3) & 9.58(-3) & 1.11(-2) & 1.31(-2) & 1.42(-2) \\
		& C($\beta=1$)     & 1.46     & 1.67     & 1.85     & 2.13     & 2.48     & 2.69     \\ \midrule
		10 & H($\beta=1/4$)   & 1.83(-3) & 2.11(-3) & 2.33(-3) & 2.71(-3) & 3.18(-3) & 3.45(-3) \\ 
		& C($\beta=1$)     & 1.47     & 1.69     & 1.87     & 2.16     & 2.53     & 2.75     \\ \midrule
		15 & H($\beta=1/4$)   & 7.75(-4) & 8.93(-4) & 9.87(-4) & 1.15(-3) & 1.34(-3) & 1.46(-3) \\
		& C($\beta=1$)     & 1.47     & 1.69     & 1.87     & 2.17     & 2.55     & 2.76     \\ \midrule
		20 & H($\beta=1/4$)   & 4.17(-4) & 4.80(-4) & 5.31(-4) & 6.16(-4) & 7.23(-4) & 7.85(-4) \\ 
		& C($\beta=1$)     & 1.47     & 1.70     & 1.88     & 2.18     & 2.55     & 2.77     \\ \bottomrule
	\end{tabular}
\end{table}

\begin{figure}[htb!]
	\centering
	\begin{minipage}[b]{0.48\textwidth}
		\includegraphics[width=\textwidth]{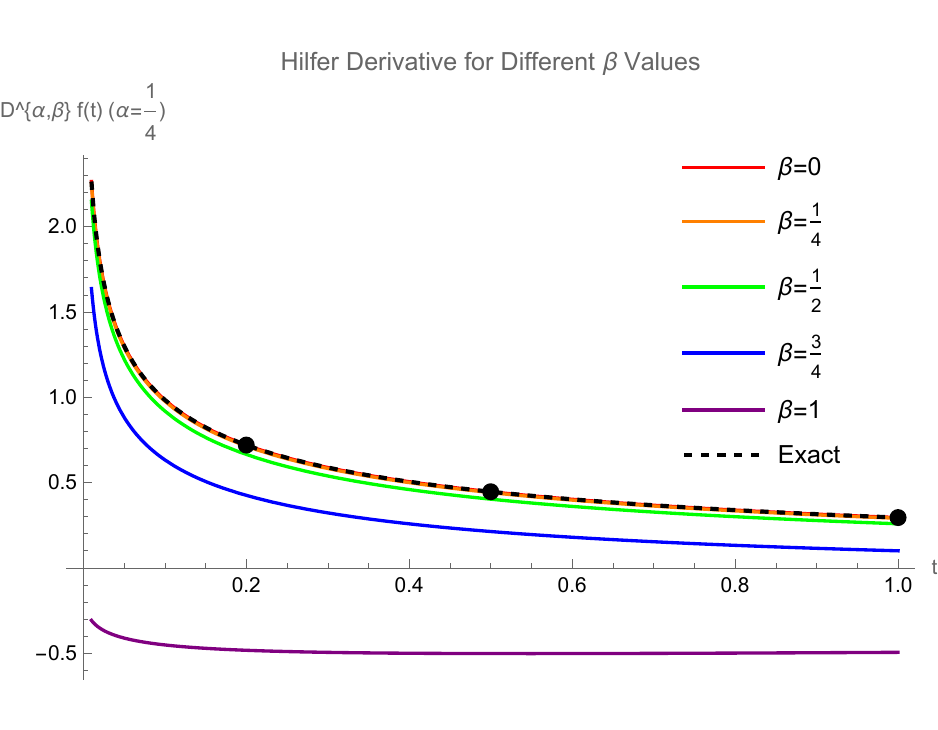}
		\caption{Exact and approximation for Ex.\ref{ex2} for different values of $\beta$ with $\alpha=1/4$, $s=1/2$ and $n=m=5$.}
		\label{fig.7}
	\end{minipage}
	\hfill
	\begin{minipage}[b]{0.48\textwidth} 
		\includegraphics[width=\textwidth]{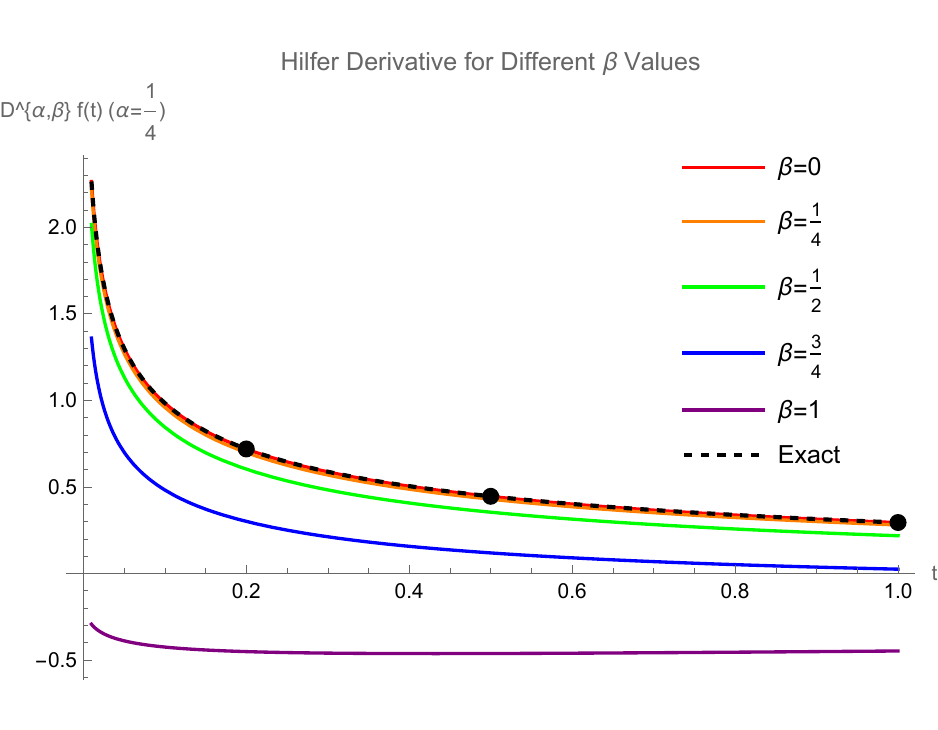}
		\caption{Exact and approximation for Ex.\ref{ex2} for different values of $\beta$ with $\alpha=1/4$, $s=1$ and $n=m=5$.}
		\label{fig.8}
	\end{minipage}
\end{figure}

\begin{table}[h]
	\centering
	\small
	\caption{Relative errors in composite Gaussian quadrature sum obtained $n$-point for fractional derivative Hilfer (H($\beta=1/4$)) and Caputo(C($\beta=1$)) in some selected values of $t \in [0, 1]$ for $\alpha = 1/4$ and $s = 1/4$.}\label{table3}
	\begin{tabular}{llllllll}
		\toprule
		\textbf{$n=m$} & \textbf{FD}  &   \textbf{$t=0.1$} & \textbf{$t=0.2$} & \textbf{$t=0.3$} & \textbf{$t=0.5$} & \textbf{$t=0.8$} & \textbf{$t=1.0$} \\ \midrule
		5  & H($\beta=1/4$)   & 1.70(-3) & 1.95(-3) & 2.16(-3) & 2.51(-3) & 2.94(-3)  & 3.19(-3)  \\
		& C($\beta=1$)     & 1.47     & 1.69     & 1.87     & 2.17     & 2.54      & 2.75      \\ \midrule
		10 & H($\beta=1/4$)   & 1.42(-4) & 1.63(-4) & 1.81(-4) & 2.10(-4) & 2.46(-4)  & 2.67(-4)  \\ 
		& C($\beta=1$)     & 1.47     & 1.70     & 1.88     & 2.18     & 2.55      & 2.77      \\ \midrule
		15 & H($\beta=1/4$)   & 2.90(-5) & 3.33(-5) & 3.69(-5) & 4.28(-5) & 5.02(-5)  & 5.45(-5)  \\
		& C($\beta=1$)     & 1.47     & 1.70     & 1.88     & 2.18     & 2.56      & 2.78      \\  \midrule
		20 & H($\beta=1/4$)   & 8.98(-6) & 1.03(-5) & 1.14(-5) & 1.33(-5) & 1.56(-5)  & 1.69(-5)  \\ 
		& C($\beta=1$)     & 1.48     & 1.70     & 1.88     & 2.18     & 2.56      & 2.78      \\ \bottomrule
	\end{tabular}
\end{table}
For Figure \ref{fig.7}, we can see that when taking $\alpha=1/4$ and $s=1/2$, the exact and approximation for different values of $\beta$ have a good matching for the Riemann derivative ($\beta=0$) and the Hilfer derivative for some values of $\beta$, but not matching for the Caputo, and need more points.

When we choose $s=1$, with five points, Figure \ref{fig.8} shows that the approximation is very close, which is exact for the Riemann derivative ($\beta=0$). However, this happens for the Hilfer derivative only for the case $\beta=1/4$.
\begin{example} \label{ex3}
	Let us consider the function $f(t)=sin(\sqrt{t})$.Then by applying the Hilfer fractional derivative, we obtain the following solution:
		\[	D^{\alpha,\beta}_{0+}  \sin\sqrt{t}
	= \frac{\sqrt{\pi}}{2} (4t)^{\frac14 - \frac{\alpha}{2}} \, J_{\frac12 - \alpha}(\sqrt{t}),
	\]
	where $J_{\alpha}(t)$ is Bessel function. 
\end{example}
In Table \ref{table4} for $\alpha=1/4$ and $s=1/2$ for some values from the interval $[0,1]$, we see that Caputo converges slowly, but for Hilfer, when taking $\beta=1/4$, we can see that convergence is satisfactory, for the Riemann derivative for the same choices in Figure \ref{fig9} the convergence is fast. 

We notice that when choosing $\beta=1/3$ for Hilfer and the same values of $\alpha$ and $s$, it shows almost the same errors with the Riemann derivative, and it converges fast.
\begin{table}[h]
	\centering
	\small
	\caption{Relative errors in quadrature sums obtained by $n$-point Riemann and Caputo rules in some selected values of $t \in [0, 1]$ for $\alpha = 1/4$ and $s = 1/2$.} \label{table4}
	\begin{tabular}{llllllll}
		\toprule
		\textbf{$n$} & \textbf{FD}  &   \textbf{$t=0.1$} & \textbf{$t=0.2$} & \textbf{$t=0.3$} & \textbf{$t=0.5$} & \textbf{$t=0.8$} & \textbf{$t=1.0$} \\ \midrule
		5  & H($\beta=1/4$)   & 1.53(-5) & 1.57(-5) & 1.60(-5) & 1.67(-5) & 1.78(-5) & 1.85(-5) \\
		& C($\beta=1$)     & 2.43(-2) & 2.48(-2) & 2.53(-2) & 2.63(-2) & 2.81(-2) & 2.93(-2) \\ \midrule
		10 & H($\beta=1/4$)   & 1.05(-6) & 1.07(-6) & 1.09(-6) & 1.14(-6) & 1.21(-6) & 1.27(-6) \\ 
		& C($\beta=1$)     & 7.13(-3) & 7.28(-3) & 7.43(-3) & 7.74(-3) & 8.25(-3) & 8.62(-3) \\ \midrule
		15 & H($\beta=1/4$)   & 2.06(-7) & 2.10(-7) & 2.15(-7) & 2.24(-7) & 2.38(-7) & 2.49(-7) \\
		& C($\beta=1$)     & 3.36(-3) & 3.43(-3) & 3.50(-3) & 3.65(-3) & 3.89(-3) & 4.06(-3) \\ \midrule
		20 & H($\beta=1/4$)   & 6.37(-8) & 6.50(-8) & 6.64(-8) & 6.92(-8) & 7.38(-8) & 7.70(-8) \\ 
		& C($\beta=1$)     & 1.94(-3) & 1.98(-3) & 2.03(-3) & 2.11(-3) & 2.25(-3) & 2.35(-3) \\ \bottomrule
	\end{tabular}
\end{table}

Figure \ref{fig11} shows the exact and approximate solutions for Example \ref{ex3} at $s=1/2$ and $n=m=5$ points, for different values of $\beta$. For $\beta=1$, where it pertains to the Caputo derivative, the estimated values do not coincide with the exact solution. For all other values, it coincides with the exact solution.Figure \ref{fig12}, with $s=1$, shows that some values of $\beta$ do not coincide the exact solution.
\begingroup
\captionsetup{skip=0pt}

\begin{figure}[htb!]
	\centering
	\begin{minipage}[b]{0.48\textwidth}
		\includegraphics[width=\textwidth, trim=0 45 0 10, clip]{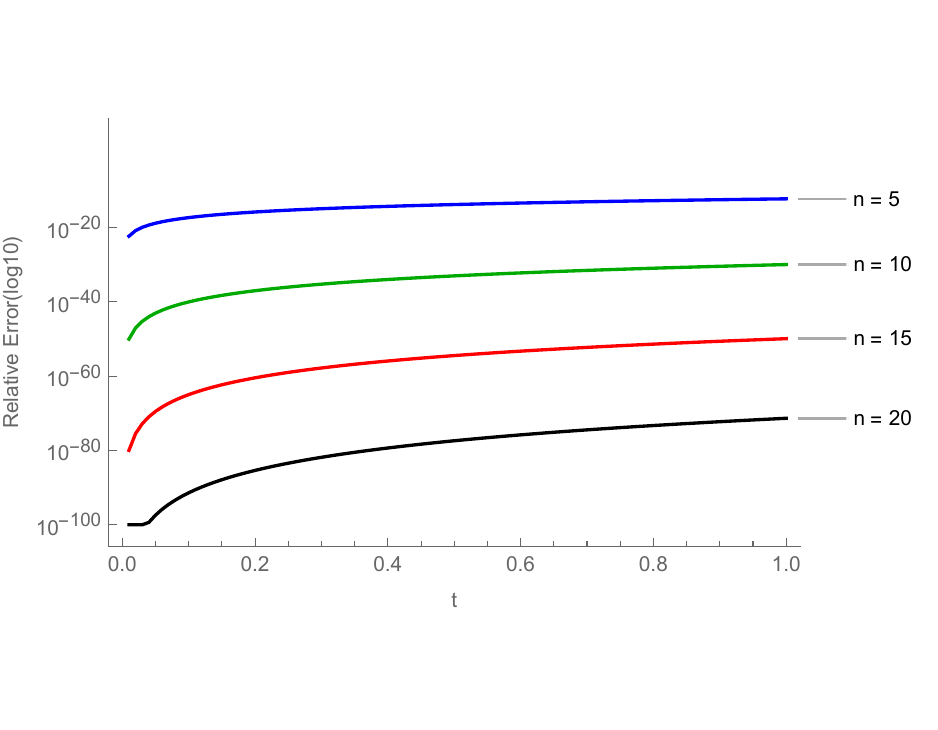}
			\caption{relative errors for Riemann derivative ($\beta=0$) for Ex. \ref{ex3} when  $\alpha=1/4$, $s=1/2$ and different values of n.}
		\label{fig9}
	\end{minipage}
	\hfill
	\begin{minipage}[b]{0.48\textwidth}
		\includegraphics[width=\textwidth, trim=0 10 0 10, clip]{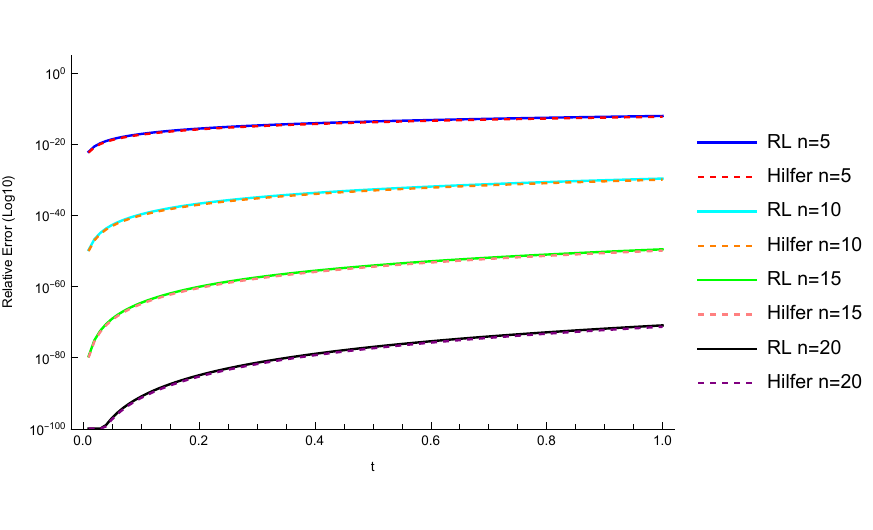}
		\caption{relative errors for ($\beta=0$,$\beta=1/3$) for Ex. \ref{ex3} when  $\alpha=1/4$, $s=1/2$ and different values of n.}
		\label{fig10}
	\end{minipage}

\end{figure}
\endgroup
\begin{figure}[htb!]
	\centering
	\begin{minipage}[b]{0.48\textwidth}
		\includegraphics[width=\textwidth]{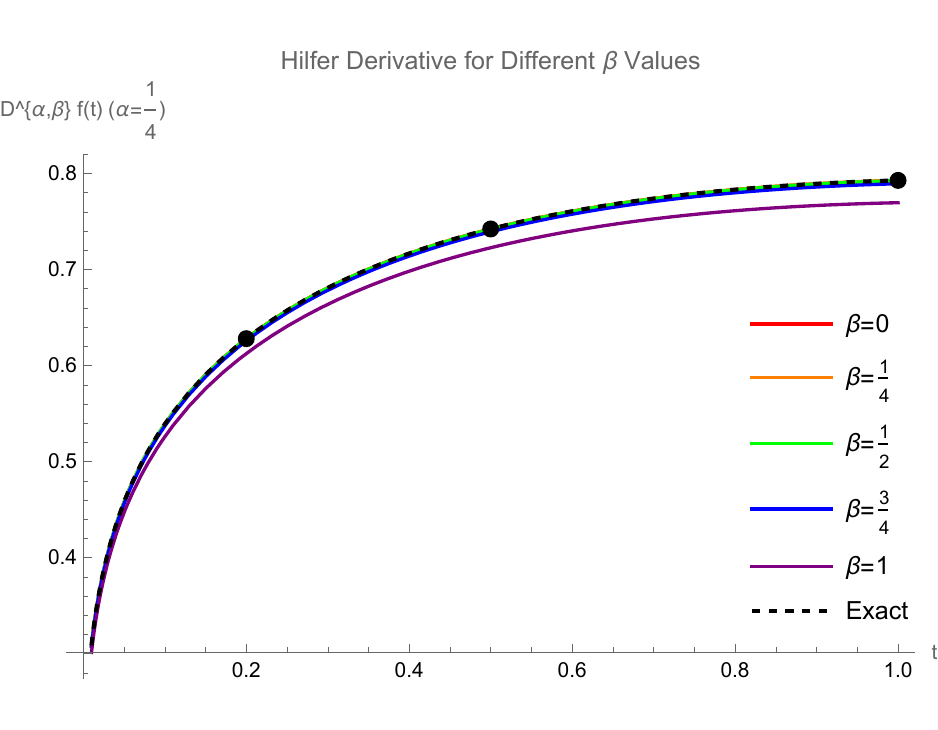}
		\caption{Exact and approximation for Ex.\ref{ex3} for different values of $\beta$ with $\alpha=1/4$, $s=1/2$ and $n=m=5$.}
		\label{fig11}
	\end{minipage}
	\hfill
	\begin{minipage}[b]{0.48\textwidth} 
		\includegraphics[width=\textwidth]{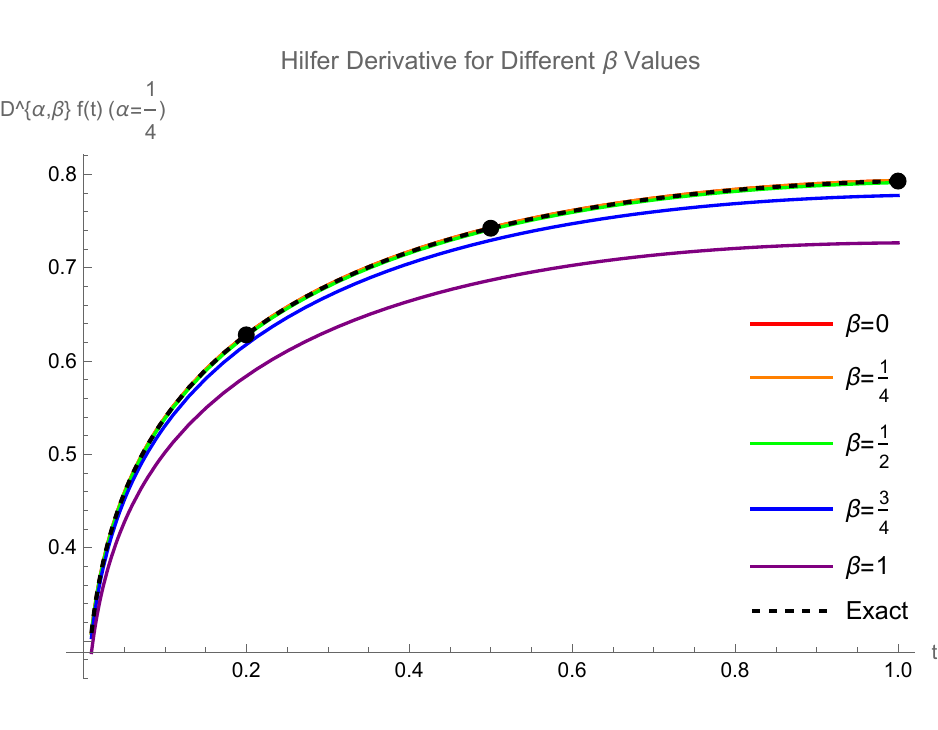}
		\caption{Exact and approximation for Ex.\ref{ex3} for different values of $\beta$ with $\alpha=1/4$, $s=1$ and $n=m=5$.}
		\label{fig12}
	\end{minipage}
\end{figure}
\section{conclusion}
In This work, we studied a numerical method for approximating the 
Hilfer fractional derivative using a composite Gaussian quadrature sum. The existing techniques have been extended to the Hilfer operator, which generalizes both the Riemann--Liouville and Caputo derivatives, and we have shown that how this valuable method applies to the fractional derivative, especially for the Riemann fractional derivative. Numerical tests confirm the method's accuracy and efficiency, as well as control of the convergence of the numerical processes by the parameters $s$ and $\beta$ in the weight function. It appears that no other method exists in the literature to approximate the Hilfer derivative operator of a function to be compared 
with the present method. While methods by \cite{S. Esmaeili_2014_Nonstandard Gauss-Lobatto quadrature,C.Li_2012_Spectral approximations to the fractional integral} are used for Riemann-Liouville and Caputo fractional derivatives, they frequently have trouble with singular functions or are computationally costly. More recently, \cite{H.KhosravianArab_2024_The sine and cosine diffusive representations} proposed diffusive representations for the Caputo fractional derivative that perform well for smooth functions, but do not address singular functions and are restricted to the Caputo case only. Our process uses fewer nodes and weights to achieve greater precision for singular cases, providing a more effective unified scheme for the Hilfer derivative. We also note that for the Caputo case $\beta = 1$, the present method performs well when the derivative $f'$ is smooth at the origin, as confirmed in Example~\ref{ex1}, but converges slowly for functions with singularities such as examples~\ref{ex2} and~\ref{ex3}, where $f'$ has a singularity at $t = 0$.

\end{document}